\documentclass[twoside,11pt]{amsart}

\usepackage{amsmath,amssymb,amsbsy,amstext,amsxtra,latexsym}

\newtheorem{theorem}{Theorem}[section]
\newtheorem{lemma}{Lemma}[section]
\newtheorem{proposition}{Proposition}[section]

\newcommand{\mF}{\mathbb F}
\def\mP{{\mathbb P}}
\def\mQ{{\mathbb Q}}
\def\mC{{\mathbb C}}

\def\tX{{\tilde X}}
\def\mZ{{\mathbb Z}}
\def\cO{{\mathcal O}}
\def\cX{{\mathcal X}}
\def\cY{{\mathcal Y}}

\def\cV{{\mathcal V}}

\def\n{\noindent}

\def\m{\medskip}

\def\t{\tilde}
\def\tX{\tilde X}
\def\T{{\mathbb T}}
\def\Ext{\operatorname{Ext}}

\def\log{\operatorname{log}}

\begin{document}

\title[Complex structure on the rational blowdown of sections in $E(4)$]
{Complex structure on the rational \\ blowdown of sections in $E(4)$}

\author{Yongnam Lee}

\address{Department of Mathematics, Sogang University,
         Sinsu-dong, Mapo-gu, Seoul 121-742, Korea}

\email{ynlee@sogang.ac.kr}

\date{Mar 11, 2010}



\begin{abstract}
We show that there is a complex structure on the symplectic
4-manifold $W_{4, k}$ obtained from the elliptic surface E(4) by
rationally blowing down $k$ sections for $2\le k\le 9$. And we
interpret it via $\mQ$-Gorenstein smoothing. This answers
affirmatively to a question raised by R. Gompf.
\end{abstract}

\maketitle

\section{Introduction}
\label{sec-1}

Let $E(n)$ denote a simply connected, relatively minimal elliptic
surface with topological Euler characteristic $c_2=12n>0$ and with a
section. The diffeomorphism type of $E(n)$ is unique, and $E(n)$ is
symplectically isomorphic to the fiber sum of $n$ copies of a
rational elliptic surface $E(1)$. Recall that $E(1)$ is obtained
from $\mP^2$ by blowing up the nine base points of a pencil of
cubics. It is known that $E(n)$ admits nine $(-n)$-curves as
disjoint sections.

Consider the case $E(4)$, and rationally blow down $k$ disjoint
sections to obtain symplectic 4-manifolds $W_{4,k}$. Rational
blowdown of $k$ disjoint sections in $E(4)$ is the same as the
normal connected sum of $E(4)$ with $k$ copies of $\mP^2$
identifying a conic in each $\mP^2$ with a section in $E(4)$. The
manifold $W_{4, 1}$ does not admit any complex structure because it
violates the Noether inequality, $p_g\le \frac{1}{2}K^2+2$ for a
minimal surface of general type: $p_g(W_{4,1})=3$ and
$K^2_{W_{4,1}}=1$.

In the paper \cite{Gom}, Gompf constructed a family of symplectic
$4$-manifolds by taking a fiber sum of other symplectic
$4$-manifolds, and he raised the following question.

\m

\n{\bf Question.} Is it possible to give a complex structure on
$W_{4, k}$ for $2\le k\le 9$?

\m

The case of $k=2$ was treated by using $\mQ$-Gorenstein smoothing
theory \cite{LP3}. We denote the singular projective surface
obtained by contracting $k$ disjoint sections in $E(4)$ by
$W_{4,k}'$. It is known that if there is a $\mQ$-Gorenstein
smoothing of $W_{4,k}'$, then a general fiber is an algebraic
surface which is symplectically isomorphic to $W_{4, k}$. Recall the
idea of the proof for the case of $k=2$ in \cite{LP3}: We consider
$E(4)$ as a double cover of Hirzebruch surface $\mF_4$ branched over
an irreducible nonsingular curve $D$ in the linear system
$|4(C_0+4f)|$, where $C_0$ is the negative section and $f$ is a
fiber of $\mF_4$. Since $D$ does not intersect $C_0$, $W_{4,2}'$ is
a double cover of the cone $\hat \mF_4$ obtained by contracting
$C_0$ in $\mF_4$. Note that $\hat \mF_4$ has a $\mQ$-Gorenstein
smoothing whose general fiber is $\mP^2$. It is obtained by a pencil
of hyperplane section of the cone of the Veronese surface imbedded
in $\mP^5$. Then $W_{4, 2}'$ has a $\mQ$-Gorenstein smoothing that
is compatible with the $\mQ$-Gorenstein smoothing of $\hat \mF_4$.
And the double covering structure extends to the $\mQ$-Gorenstein
smoothing.

\m

In this paper, we will answer affirmatively to the above question.
We will construct directly a $\mQ$-Gorenstein smoothing of $W_{4,
k}'$ for $2\le k\le 9$.

\m \n{\bf Theorem.} It is possible to give a complex structure on
the rational blowdown of $k$ sections in $E(4)$ for $2\le k\le 9$.

\m

The construction is as follows: A bidouble cover of $E(1)$ branched
over three general fibers becomes a $E(4)$ with nine disjoint
sections. Choose two smooth cubics $D_1, D_2$ in $\mP^2$ such that
$D_1$ and $D_2$ meet transversally at nine points $p_1, \ldots,
p_9$. Let $D_0$ be a general member of the pencil induced by $D_1$
and $D_2$. Let $V$ be a bidouble cover of $\mP^2$ branched over
$D_1, D_2$, and $D_0$. Then $V$ has nine $1/4(1,1)$ singularities,
and its minimal resolution is $E(4)$.  By moving the cubic $D_0$,
one can construct a one-parameter family of smooth cubics $D_t$ such
that $D_t$ does not through $p_1, \ldots, p_k$ and $D_t$ intersects
$D_1, D_2$ transversally at $p_{k+1}, \ldots, p_9$ if $t \ne 0$.
First, we note that $k\ge 2$, because if a cubic passes through 8
intersection points then it also passes through the 9-th point. This
also explains that $W_{4, 1}$ does not admit any complex structure.
One-parameter family of smooth cubics $D_t$ induces a one-parameter
family $\cV_k$ by using a bidouble cover of $\mP^2$ branched over
$D_1, D_2, D_t$. Then by using a simultaneous resolution, we have a
one-parameter family $\cX$ such that $X_0=W_{4, k}'$ and a general
fiber $X_t$ is smooth with $K_{X_t}^2=k$. This family is a
$\mQ$-Gorenstein smoothing of $W_{4, k}'$.

In Section 3, we briefly review the theory of $\mQ$-Gorenstein
smoothing and interpret the result via the local-global exact
sequence of smoothings.

\m

{\em Acknowledgements}. The author would like to thank Jongil Park
for information about the question raised by Gompf, and to thank
JongHae Keum and Miles Reid for suggesting some comments to improve
the paper. The author was supported by the WCU Grant funded by the
Korean Government (R33-2008-000-10101-0). The part of the paper was
worked out during his visit to MSRI in January, 2009. He thanks MSRI
for hospitality during his visit and for supporting travel expenses.

\section{Proof of Theorem}
\label{sec-2}

A bidouble cover is a $\mZ_2\oplus\mZ_2$ cover $\psi: V \to U$, with
$U$ a smooth surface. The building data for a bidouble cover consist
of the following:
\begin{enumerate}
\item smooth divisors $D_1, D_2, D_3$ in $U$ having pairwise
transverse intersections, and
\item line bundles $L_1, L_2, L_3$ such that $2L_g= D_j+D_k$ for
each permutation $(g, j, k)$ of $(1,2,3)$.
\end{enumerate}

We have a decomposition
\[\psi_*\cO_V=\cO_U\oplus \sum_{\chi\in G^*} L_{\chi}^{-1},\]
where $L_\chi$ is a line bundle and $G=\mZ_2\oplus\mZ_2$ acts on
$L_\chi^{-1}$ via the character $\chi\in G^*$. Since we have three
elements $g_1, g_2, g_3$ in $G^*$, there are three characters
$\chi_i\in G^*$ for $i=1,2,3$.

Suppose that $D_1, D_2$, and $D_3$ have no common intersection. Then
$V$ is also smooth, and
\[ p_g(V)=p_g(U)+\sum_{i=1}^3 h^0(U, \cO_U(K_U+L_{\chi_i})),\]
\[ \chi(\cO_V)=4\chi(\cO_U)+\frac{1}{2}\sum_{i=1}^3L_{\chi_i}(K_U+L_{\chi_i}),\]
\[ 2K_V=\psi^*(2K_U+L_{\chi_1}+L_{\chi_2}+L_{\chi_3}).\]
For more information on bidouble covers, see \cite{Ca}, or
\cite{Man94}, or \cite{Pa}.

\begin{lemma}
\label{lem-2.1} If $D_1, D_2$, and $D_3$ have a common intersection
point $p$ in $X$ then $V$ has the singularity $1/4(1,1)$.
\end{lemma}

\begin{proof}
Since the problem is local, we may assume that $U=\mC^2$, $D_1:
x_1=0$, $D_2: x_2=0$, and $D_3: x_1+x_2=0$ where $x_1, x_2$ are
coordinates of $\mC^2$. Let
\[ y_1^2=x_1, \, y_2^2=x_2, \, y_3^2=x_1+x_2.\]
Since $y_1^2+y_2^2=y_3^2$, the simple cover of $\mC^2$ induced by
$\mZ_2\oplus\mZ_2\oplus\mZ_2$ has a $A_1$-singularity. And a
bidouble cover is the quotient of the above simple cover by the
involution multiplying all the three functions by -1. Therefore $V$
has a singularity $1/4(1,1)$.
\end{proof}

Let $D_1, D_2$ be two smooth cubics in $\mP^2$ such that $D_1$ and
$D_2$ meet transversally at nine points $p_1, \ldots, p_9$. Let
$D_0$ be a general member of the pencil induced by $D_1$ and $D_2$.
Let $V$ be a bidouble cover of $\mP^2$ branched over $D_1, D_2$, and
$D_0$. Then $V$ has nine $1/4(1,1)$ quotient singularities by Lemma
2.1, and $p_g(V)=3, \chi(\cO_V)=4, K_V^2=9$. Its minimal resolution
is $E(4)$.

By moving the cubic $D_0$, one can construct a one-parameter family
of smooth cubics $D_t$ such that $D_t$ does not through $p_1,
\ldots, p_k$ and $D_t$ intersects $D_1, D_2$ transversally at
$p_{k+1}, \ldots, p_9$ if $t \ne 0$. First, we note that $k\ge 2$,
because if a cubic passes through 8 intersection points then it also
passes through the 9-th point. One-parameter family of smooth cubics
$D_t$ induces a one-parameter family $\cV_k$ by using a bidouble
cover of $\mP^2$ branched over $D_1, D_2, D_t$. Recall the
definition of a $\mQ$-Gorenstein smoothing (cf. \cite{KSB}).

\m

\n{\bf Definition.} Let $X$ be a normal projective surface with
quotient singularities. Let $\cX\to\Delta$ (or $\cX/\Delta$) be a
flat family of projective surfaces over a small disk $\Delta$. The
one-parameter family of surfaces $\cX\to\Delta$ is called a {\it
$\mQ$-Gorenstein smoothing} of $X$ if it satisfies the following
three conditions;

\n(i) the general fiber $X_t$ is a smooth projective surface,

\n(ii) the central fiber $X_0$ is $X$,

\n(iii) the canonical divisor $K_{\cX/\Delta}$ is $\mQ$-Cartier.

\m

By using a simple cover $\mZ_2\oplus\mZ_2\oplus\mZ_2$ and by
simultaneous resolution of $A_1$-singularities, we have a
one-parameter family $\cX$ over a small disk $\Delta$ in $\mC$ such
that $X_0=W_{4, k}'$ and a general fiber $X_t$ is smooth with
$K_{X_t}^2=k$. This family is a $\mQ$-Gorenstein smoothing of $W_{4,
k}'$ because $2K_{\cX/\Delta}$ is Cartier. It is well-known that a
general fiber of a $\mQ$-Gorenstein smoothing is the rational
blowdown of $k$ disjoint sections in $E(4)$. This gives the proof of
Theorem.

\m \n{Remark.} One can also do by using a bidouble cover of
$\mP^1\times\mP^1$ branched over $D_1, D_2$, and $D_3$ in
$H^0(\mP^1\times\mP^1, \cO_{\mP^1\times\mP^1}(2, 2))$. Then one can
show that there is a complex structure on the rational blowdown of
$k$ sections in $E(4)$ for $2\le k\le 8$. Since a cubic in $\mP^2$
can be transformed to a (2, 2) divisor in $\mP^1\times\mP^1$ by
blowing up twice and blowing down once, two constructions of
$\mQ$-Gorenstein smoothing of $W_{4, k}'$ can be identified.

\section{$\mQ$-Gorenstein smoothing}
\label{sec-3}

In this section we briefly review the theory of $\mQ$-Gorenstein
smoothing of projective surfaces with special quotient
singularities. The $\mQ$-Gorenstein smoothing theory can be used to
construct surfaces of general type, but the methodology comes from
3-fold Mori theory (results of Koll\'ar and Shepherd-Barron on
$\mQ$-Gorenstein smoothing and relative canonical models). The
compactification theory of a moduli space of surfaces of general
type was established during the last 20 years. It was originally
suggested in \cite{KSB} and it was established by Alexeev's proof
for boundness~\cite{Alex} and by the Mori program for threefolds
(\cite{KM} for details). It is natural to expect the existence of a
surface with special quotient singularities in the boundary of a
compact moduli space.

\m

A quotient singularity which admits a $\mQ$-Gorenstein smoothing is
called a {\it singularity of class T}.

\begin{proposition}[\cite{KSB, Man91, Wa1}]
\label{pro-3.1} Let $(X_0, 0)$ be a germ of two dimensional quotient
singularity. If $(X_0, 0)$ admits a $\mQ$-Gorenstein smoothing over
the disk, then $(X_0, 0)$ is either a rational double point or a
cyclic quotient singularity of type $\frac{1}{dn^2}(1, dna-1)$ for
some integers $a, n, d$ with $a$ and $n$ relatively prime.
\end{proposition}

\begin{proposition}[\cite{KSB, Man91, Wa2}]
\label{pro-3.2}
\begin{enumerate}
\item The quotient singularities ${\overset{-4}{\circ}}$ and
${\overset{-3}{\circ}}-{\overset{-2}{\circ}}-{\overset{-2}{\circ}}-\cdots-
{\overset{-2}{\circ}}-{\overset{-3}{\circ}}$ are of class $T$.
\item If the singularity
${\overset{-b_1}{\circ}}-\cdots-{\overset{-b_r}{\circ}}$ is of class
$T$, then so are
$${\overset{-2}{\circ}}-{\overset{-b_1}{\circ}}-\cdots-{\overset{-b_{r-1}}
{\circ}}- {\overset{-b_r-1}{\circ}} \quad\text{and}\quad
{\overset{-b_1-1}{\circ}}-{\overset{-b_2}{\circ}}-\cdots-
{\overset{-b_r}{\circ}}-{\overset{-2}{\circ}}.$$
\item Every singularity of class $T$ that is not a rational double point can be
obtained by starting with one of the singularities described in
$(1)$ and iterating the steps described in $(2)$.
\end{enumerate}
\end{proposition}

Let $X$ be a normal projective surface with singularities of class
$T$. Due to the result of Koll\'ar and Shepherd-Barron \cite{KSB},
there is a $\mQ$-Gorenstein smoothing locally for each singularity
of class $T$ on $X$ (see Proposition 3.5). The natural question is
whether this local $\mQ$-Gorenstein smoothing can be extended over
the global surface $X$ or not. The answer can be obtained by
figuring out the obstruction map of the sheaves of deformation
$T^i_X=Ext^i_X(\Omega_X,\cO_X)$ for $i=0,1,2$. For example, if $X$
is a smooth surface, then $T^0_X$ is the usual holomorphic tangent
sheaf $T_X$ and $T^1_X=T^2_X=0$. By applying the standard result of
deformations \cite{LS, Pal} to a normal projective surface with
quotient singularities, we get the following

\begin{proposition}[\cite{Wa1}, \S 4]
\label{pro-3.3} Let $X$ be a normal projective surface with quotient
singularities. Then
\begin{enumerate}
\item The first order deformation space of $X$ is represented by
the global Ext 1-group $\T^1_X=\Ext^1_X(\Omega_X, \cO_X)$.
\item The
obstruction space lies in the global Ext 2-group
$\T^2_X=\Ext^2_X(\Omega_X,\cO_X)$.
\end{enumerate}
\end{proposition}

Furthermore, by applying the general result of local-global spectral
sequence of ext sheaves (\cite{Pal}, \S 3) to deformation theory of
surfaces with quotient singularities so that $E_2^{p, q}=H^p(T^q_X)
\Rightarrow \T^{p+q}_X$, and by $H^j(T^i_X)=0$ for $i, j\ge 1$, we
also get

\begin{proposition}[\cite{Man91, Wa1}]
\label{pro-3.4} Let $X$ be a normal projective surface with quotient
singularities. Then
\begin{enumerate}
\item We have the exact sequence
$$0\to H^1(T^0_X)\to \T^1_X\to \ker [H^0(T^1_X)\to H^2(T^0_X)]\to 0$$
where $H^1(T^0_X)$ represents the first order deformations of $X$
for which the singularities remain locally a product.
\item If $H^2(T^0_X)=0$, every local deformation of
the singularities may be globalized.
\end{enumerate}
\end{proposition}

The vanishing $H^2(T^0_X)=0$ can be obtained via the vanishing of
$H^2(T_V(-\log \ E))$, where $V$ is the minimal resolution of $X$
and $E$ is the reduced exceptional divisors. Note that every
singularity of class $T$ has a local $\mQ$-Gorenstein smoothing by
Proposition 3.5 below. With the help of the birational geometry of
threefolds and their applications to deformations of surface
singularities, the following proposition is obtained.
 Note that the cohomology $H^0(\t{T}^1_X)$ is given explicitly.

\begin{proposition}[\cite{KSB, Man91}]
\label{pro-3.5}
\begin{enumerate}
\item Let $a, d, n>0$ be integers with $a, n$ relatively prime and
consider a map $\pi :\cY/\mu_n \to \mC^d$, where $\cY\subset
\mC^3\times \mC^d$ is the hypersurface of equation
$uv-y^{dn}=\sum_{k=0}^{d-1}t_ky^{kn};\ t_0,\ldots, t_{d-1}$ are
linear coordinates over $\mC^d$, $\mu_n$ acts on $\cY$ by
$$\mu_n\ni\xi : (u, v, y, t_0,\ldots, t_{d-1})\to (\xi u, \xi^{-1}v,
\xi^ay, t_0, \ldots, t_{d-1})$$ and $\pi$ is the factorization
through the quotient of the projection $\cY\to\mC^d$. Then $\pi$ is
a $\mQ$-Gorenstein smoothing of the cyclic singularity germ $(X_0,
0)$ of type $\frac{1}{dn^2}(1, dna-1)$. Moreover every
$\mQ$-Gorenstein smoothing of $(X_0, 0)$ is isomorphic to the
pull-back of $\pi$ under some germ of holomorphic map $(\mC, 0)\to
(\mC^d, 0)$.
\item Let $X$ be a normal projective surface with singularities of class T.
Then
$$H^0(\t{T}^1_X)=\sum_{p\,\in\text{\,singular points of $X$}}
\mC_p^{\oplus d_p}$$ where $p$ is a point of type $\frac{1}{d_p n^2}
(1, d_p an-1)$ with $(a, n)=1$.
\end{enumerate}
\end{proposition}

\begin{theorem}[\cite{LP1}]
\label{thm-3.1} Let $X$ be a normal projective surface with
singularities of class $T$. Let $\pi: V\to X$ be the minimal
resolution and let $E$ be the reduced exceptional divisors. Suppose
that $H^2(T_V(-\log \ E))=0$. Then $H^2(T^0_X)=0$ and there is a
$\mQ$-Gorenstein smoothing of $X$.
\end{theorem}

As we see in Theorem 3.1 above, if $H^2(T^0_X)=0$, then there is a
$\mQ$-Gorenstein smoothing of $X$. For example, a simply connected
minimal surface of general type with $p_g=0$ and $K^2 =2$ was
constructed in \cite{LP1} by proving $H^2(T^0_X)=0$. The vanishing
$H^2(T^0_X)=0$ is obtained by the careful study of the configuration
of singular fibers in a rational elliptic surface.

\m

\n{\bf Remark.} All recent constructions of surfaces of general type
with vanishing geometric genus via $\mQ$-Gorenstein smoothing \cite
{LP1} \cite{LP2} \cite{PPS1} \cite{PPS2} \cite{PPS3} have the
vanishing $H^2(T^0_X)=0$ by the same arguments in \cite{LP1}. By
upper semi-continuity, $H^2(X_t, T_{X_t})=0$ for a general fiber
$X_t$ of a $\mQ$-Gorenstein smoothing. We note that
\[h^1(X_t, T_{X_t})-h^2(X_t, T_{X_t})=10\chi(\cO_{X_t})-2K_{X_t}^2.\]
If $H^2(X_t, T_{X_t})=0$ then the dimension of the deformation space
of $X_t$ is $h^1(X_t, T_{X_t})=10-2K_{X_t}^2.$ It implies that there
is no nontrivial deformation of $X_t$ if $K_{X_t}^2\ge 5$. In
particular, it is not possible to construct a minimal surface of
general type $X_t$ with $p_g=0, K^2\ge 5$ via $\mQ$-Gorenstein
smoothing of a singular surface $X$ with $H^2(T^0_X)=0$.

\m

But, in general, the cohomology $H^2(T^0_X)$ is not zero and it is a
very difficult problem to determine whether there exists a
$\mQ$-Gorenstein smoothing of $X$. Hence, in the case that
$H^2(T^0_X) \not =0$, we need to develop another technique in order
to investigate the existence of  $\mQ$-Gorenstein smoothing. Even
though we do not know whether such a technique exists in general, if
$X$ is a normal projective surface with singularities of class $T$
which admits a cyclic group action with some additional properties,
then we are able to show that it admits a $\mQ$-Gorenstein
smoothing. Explicitly, we get the following theorem.

\begin{theorem}[\cite{LP3}]
\label{thm-3.2} Let $X$ be a normal projective surface with
singularities of class $T$. Assume that a cyclic group $G$ acts on
$X$ such that
\begin{enumerate}
\item $Y=X/G$ is a normal projective
surface with singularities of $T$,
\item  $p_g(Y)=q(Y)=0$,
\item $Y$ has a $\mQ$-Gorenstein smoothing,
\item the map $\sigma: X\to Y$ induced by a cyclic covering is flat,
and the branch locus $D$ (resp. the ramification locus) of the map
$\sigma: X\ \to Y$ is nonsingular curve lying outside the singular
locus of $Y$ (resp. of $X$), and
\item $H^1(Y, \cO_Y(D))=0$.
\end{enumerate}
Then there exists a $\mQ$-Gorenstein smoothing of $X$ that is
compatible with a $\mQ$-Gorenstein smoothing of $Y$. And the cyclic
covering extends to the $\mQ$-Gorenstein smoothing.
\end{theorem}

\n {\bf Example.} Assume that $n \ge 5$ and
 let $\mF_n$ be a Hirzebruch surface. Let $C_0$ be
 a negative section with $C_0^2=-n$ and $f$ be a fiber of $\mF_n$.
 Choose a special irreducible (singular) curve
 $D$ in the linear system $|4(C_0+nf)|$ having a special intersection
 with one special fiber $f$: Note that $D\cdot f=4$.
 We want $D$ to intersect with $f$ at two distinct points $p$ and
 $q$ that are not in $C_0$.
 Let $x=0$ be the local equation of $f$ and $x, y$ be a coordinate
 at $p$ (resp. at $q$). We require that the local equation of $D$ at $p$
 (resp. at $q$) is $(y-x)(y+x)=0$ (resp. $(y-x^{n-4})(y+x^{n-4})=0$).
 In \cite{LP3}, we show that
 there is a curve $D$ satisfying the conditions above and being
 nonsingular at every point except at the two points
 $p$ and $q$.
 Let $\sigma: \tX_n\to\mF_n$ be a double cover branched over
 the curve $D$ chosen above. Then $\tX_n$ is a singular elliptic surface
 with $p_g=n-1$ and $\chi(\cO_{\tX_n}) =n$
 which has two rational double points by the local equations of $D$
 at $p$ and $q$; one is of type $A_1$ and the other one
 is of type $A_{2n-9}$.
 Therefore its minimal resolution is also an elliptic surface $E(n)$.
 First we blow up at $p$ and $q$ in $\mF_n$. Then we have an
 exceptional curve obtained by a blowing up at $p$ which intersects
 with the proper transform of $D$ transversally at two points,
 and we also have an exceptional curve obtained by a blowing up at $q$
 which intersects with the proper transform of $D$ at one point, say $q_1$.
 Let $x=0$ be the local equation of the $(-1)$-exceptional curve at
 $q_1$. Then the local equation of the proper transform of $D$ at
 $q_1$ is $(y-x^{n-5})(y+x^{n-5})=0$. We blow up again at $q_1$. By
 the continuation of blowing up at infinitely near points of $q$,
 we have the following configuration of smooth rational curves

 \[\begin{array}{cccc}
   \underset{U_{n-3}}{\overset{-n}{\circ}} - &
   \underset{U_{n-4}}{\overset{-2}{\circ}} &
   - \underset{U_{n-5}}{\overset{-2}{\circ}}-
   \cdots - &
  \underset{U_1}{\overset{-2}{\circ}}\\
   & \vert & & \vert \\
   & {\underset{E_1, -1}{\circ}} & & {\underset{E_2, -1}{\circ}}
  \end{array}\]
 where the proper transform of $D$ intersects with $E_i$, $i=1,2$
 at two points transversally. We denote this surface by $Z_n$ obtained
 by $(n-3)$ times blowing-ups of $\mF_n$. Next, by Artin's
criterion for contractibility~\cite{Art},
 we can contract a configuration $C_{n-2}$,
 which is a linear chain of $\mP^1$'s
 \[\underset{U_{n-3}}{\overset{-n}{\circ}}-\underset{U_{n-4}}
   {\overset{-2}{\circ}} -\underset{U_{n-5}}{\overset{-2}{\circ}}-
   \cdots - \underset{U_1}{\overset{-2}{\circ}} , \]
 so that it produces a singular normal projective surface.
 We denote this surface by $Y_n$.
 We note that $\Delta$ is the proper transform of $D$ in $Z_n$ and
 that $Y_n$ has a cyclic quotient singularity of type
 $\frac{1}{(n-2)^2}(1, n-3)$, which is a singularity of class $T$.
 In \cite{LP3}, it is proved that the singular surface $Y_n$ admits a $\mQ$-Gorenstein smoothing.
 Note that $\tX_n$ is a double covering of $\mF_n$ branched over $D$,
 and the minimal resolution of two rational double points of type
 $A_1$ and $A_{2n-9}$ in $\tX_n$ is $E(n)$, which is also a double
 cover of $Z_n$ branched over the proper transform of $D$.
 Since the proper transform of $D$ does not meet the contracted linear
 chain of $\mP^1$'s, we have a double cover of $Y_n$ branched
 over the image of the proper transform of $D$ by the map $\psi$.
 We denote this surface by $X_n$.
 Then $X_n$ is a singular surface obtained by contracting two disjoint
 configurations $C_{n-2}$ from an elliptic surface $E(n)$ and
 it has two quotient singularities of class $T$,
 both are of type $\frac{1}{(n-2)^2}(1, n-3)$.
 Therefore we have the following commutative diagram of maps

\[\begin{array}{ccccc}
   \tX_n &  \leftarrow & E(n) & \rightarrow & X_n \\
   \downarrow & & \downarrow & & \downarrow\\
   \mF_n &  \stackrel{\pi}{\leftarrow} & Z_n & \stackrel{\psi}{\rightarrow}
    & Y_n
 \end{array} \]

\n where all vertical maps are double covers.
 Then, it is shown in \cite{LP3} that the singular surface $X_n$ has a $\mQ$-Gorenstein
 smoothing of two quotient singularities simultaneously by Theorem
 3.2, and its smoothing is a Horikawa surface (cf. \cite{FS}).

\m

\n {\bf Remark.} Let $X=W_{4,k}'$ for $2\le k\le 9$ in Introduction.
Then $H^2(T^0_X)\ne 0$ because $H^2(T_{E(4)})\ne 0$. We consider the
local-global exact sequence of smoothing in Proposition 3.4
\[ H^0(T^1_X)=\oplus_{i=1}^k\mC_i\to H^2(T^0_X) \]
where $\mC_i=\mC$. Then the dimension of the image of the above
obstruction map is at most one because there is a $\mQ$-Gorenstein
smoothing of any $s\ge 2$ quotient singularities simultaneously by
the argument in Section 2.


\begin{thebibliography}{999}

\bibitem[1]{Alex} V. Alexeev, {\it Boundness and $K^2$ for log surfaces}, Internat. J. Math.
                {\bf 5} (1994), 779--810.

\bibitem[2]{Art} M. Artin, {\it Some numerical criteria for contractibility of
        curves on algebraic surfaces}, Amer. J. Math. {\bf 84} (1962), 485--496.

\bibitem[3]{Ca} F. Catanese, {\it Singular bidouble covers and the construction of interesting algebraic surfaces},
            Algebraic geometry: Hirzebruch 70 (Warsaw, 1998), 97--120, Contemp. Math., {\bf 241},
            Amer. Math. Soc., Providence, RI, 1999.

\bibitem[4]{FS} R. Fintushel and R. Stern, {\it Rational blowdowns of smooth
                 4-manifolds}, Jour. Diff. Geom. {\bf 46} (1997), 181--235.

\bibitem[5]{Gom} R. Gompf, {\it A new construction of symplectic manifolds},
                  Ann. of Math. (2) {\bf 142} (1995), 527--595.

\bibitem[6]{KM} J. Koll\'ar and S. Mori, Birational geometry of algebraic
                  varieties, {\bf134} Cambridge Tracts in Mathematics, 1998.

\bibitem[7]{KSB} J. Koll\'ar and N.I. Shepherd-Barron, {\it Threefolds and deformations of surface singularities},
              Invent. Math. {\bf 91} (1988), 299--338.

\bibitem[8]{LP1} Y. Lee and J. Park, {\it A simply connected surface of general
                 type with $p_g=0$ and $K^2=2$}, Invent. Math.
                  {\bf 170} (2007), 483--505.

\bibitem[9]{LP2} Y. Lee and J. Park, {\it A complex surface of general type with $p_g=0, K^2=2$ and $H_1=\mZ/2\mZ$},
                 Math. Res. Lett. {\bf 16} (2009), 323--330.
.
\bibitem[10]{LP3} Y. Lee and J. Park, {\it A construction of Horikawa surface via
                Q-Gorenstein smoothings}, to appear in Math. Z.

\bibitem[11]{LS}  S. Lichtenbaum and M. Schlessinger, {\it The cotangent complex
                  of a morphism}, Trans. Amer. Math. Soc. {\bf 128} (1967) 41--70.

\bibitem[12]{Man91} M. Manetti, {\it Normal degenerations of the complex projective
                    plane}, J. Reine Angew. Math. {\bf 419} (1991), 89-118.

\bibitem[13]{Man94} M. Manetti, {\it On some components of moduli space of surfaces of general
                    type}, Compositio Mathematica, {\bf 92} (1994),
                    285--297.

\bibitem[14]{Pal} V. P. Palamodov, {\it Deformations of complex spaces},
                  Russian Math. Surveys {\bf 31:3} (1976), 129--197.

\bibitem[15]{Pa} R. Pardini, {\it Abelian covers of algebraic
varieties}, J. Reine Angew. Math. {\bf 417} (1991), 191--213.

\bibitem[16]{PPS1} H. Park, J. Park, D. Shin, {\it A simply connected surface of general type with $p_g=0$ and $K^2=3$},
              Geometry \& Topology {\bf 13} (2009), 743--767.

\bibitem[17]{PPS2} H. Park, J. Park, D. Shin, {\it A complex surface of general type with $p_g=0, K^2=3$ and
               $H_1=Z/2Z$}, arXiv:0803.1322.

\bibitem[18]{PPS3} H. Park, J. Park, D. Shin, {\it A simply connected surface of general type with
               $p_g=0$ and $K^2=4$}, Geometry \& Topology {\bf 13} (2009),
               1483--1494.

\bibitem[19]{Wa1} J. Wahl, {\it Smoothing of normal surface singularities},
                  Topology {\bf 20} (1981), 219--246.

\bibitem[20]{Wa2} J. Wahl, {\it Elliptic deformations of minimally
                  ellitic singularities}, Math. Ann. {\bf253} (1980), 241--262.


\end{thebibliography}
\end{document}